\documentstyle[12pt,amssymb]{article}

\input epsf.tex
\begin{document}

\begin{titlepage}

\begin{center}

\phantom{sdfsdfsdfsff}

\vspace{4cm}

{\Large{\bf Quantum tori, mirror symmetry and deformation theory}}

\vspace {1cm}

{\bf  Yan Soibelman}\\
{Department of Mathematics,\\
 Kansas State University,\\
Manhattan, KS 66506, USA}

\end{center}
\end{titlepage}

\newtheorem{thm}{Theorem}
\newtheorem{lmm}{Lemma}
\newtheorem{dfn}{Definition}
\newtheorem{rmk}{Remark}
\newtheorem{prp}{Proposition}
\newtheorem{conj}{Conjecture}
\newtheorem{exa}{Example}
\newtheorem{que}{Question}

\newtheorem{ack}{Acknowledgements}
\newcommand{\K}{{\bf k}}
\newcommand{\C}{\bf C}
\newcommand{\R}{\bf R}
\newcommand{\N}{\bf N}
\newcommand{\Z}{\bf Z}
\newcommand{\G}{\Gamma}
\newcommand{\A}{A_{\infty}}

\newcommand{\ihom}{\underline{\Hom}}
\newcommand{\ra}{\longrightarrow}
\newcommand{\epi}{\twoheadrightarrow}
\newcommand{\mono}{\hookrightarrow}

\section{Introduction}

\subsection{}

Mathematical models of dualities in quantum physics is a very interesting and intriguing
area of research. It became clear after the work of
Kontsevich (see [Ko1]-[Ko3]) that the ``right'' framework for general duality theorems 
is the (yet non-existing) theory of  moduli spaces of $\A$-categories.
Informally, an $\A$-category (with extra conditions imposed) models
a ``projective non-commutative space'' together with the 
formal moduli space of its deformations.
 Moduli space of such
non-commutative spaces (whatever it is) consists of many ``connected components''.
The boundary of the compactification of a component may contain ``cusps''.
A non-commutative space can degenerate into
a commutative one at a cusp. In some cases one can assign to a pair
{\it (component, cusp)} a ``dual'' pair. The corresponding $\A$-categories are
equivalent. In particular, their Hochschild cohomologies (interpreted as
tangent spaces to the moduli space of $\A$-categories) are isomorphic.
This idea of Kontsevich has received  spectacular confirmation in homological
mirror symmetry program (see [Ko1], [KoSo2]).

\subsection{}
Even if the degeneration of a non-commutative space at a cusp 
is not commutative, one hopes to extend the dualities to the boundary stratum.
The purpose of this paper
is to explain without details the idea
of non-commutative compactification in an example
of abelian varieties.
It will be
discussed at length elsewhere (see [So]).
The paper
is an extended version of a series of talks I gave at Ecole Polytechnique,
MSRI, Stanford University, Max-Planck Institut f\"ur Mathematik in Bonn,
Oberwolfach workshop on non-commutative geometry and
Moshe Flato Euroconference in Dijon in 1999-2000. 
In this paper we are going to 
treat many aspects informally, in order to explain
main ideas.

Let us start with an example
of what can be thought of as a non-commutative compactification.

\begin{exa}
The universal covering of the moduli space of elliptic curves
admits a ``non-commutative" compactification with the boundary stratum
consisting of the universal covering
of the moduli space of quantum (= non-commutative) tori.
Moreover, the $SL(2,{\bf Z})$-symmetry extends from the
upper-half plane to the boundary real line. 

\end{exa}

We understand
a ``space" as a category of certain sheaves on it. Thus 
the ``moduli space" of ``spaces" corresponds to the  ``moduli space" 
of categories. 
The latter ``moduli space" ${\cal M}$ can be
a usual manifold or orbifold  (a compactification of ${\cal M}$
is often a compact manifold  with corners). 
Let us consider a path $\gamma:[0,1]\to {\cal M}$
between an interior point of the
compactified moduli space and a point
of the boundary. Then we have a 1-parameter family of categories ${\cal C}_t$
along the path. It can happen that ${\cal C}_t, t\in [0,1)$
is a category of sheaves on a topological space, while ${\cal C}_1$
is not.
We are going to treat ${\cal C}_1$
as a category of sheaves on a ``non-commutative"
space. Such non-commutative spaces consititute a ``non-commutative stratum"
of the compactification $\overline{\cal M}$.

In the example above, we think about an elliptic curve
$E_q={\bf C}^{\ast}/q^{\bf Z}, q=e^{2\pi i\tau}, Im (\tau)>0$ as about 
the bounded derived category $D^b(E_q)$ of coherent
sheaves on it (we are going to discuss below reasons why derived categories
appear in the story).
Let us consider the {\it universal covering} of the moduli space
of complex structures.
 Thus we have a family of categories $D^b(E_q)$ parametrized
by the upper-half plane ${\cal H}=\{\tau|Im(\tau)>0\}$. 
As $Im(\tau)=0$ the elliptic curve does not exist as a ``commutative" space.
We will argue that the degenerate object is represented by the
derived category of certain modules over the algebra of functions
on a quantum torus.
From this point of view elliptic curves and quantum tori belong to the same
family of non-commutative spaces. Also, the $SL(2,{\bf Z})$-symmetry should be present
for all $q\ne 0$. See Section 3 for more details.
\footnote{The idea to treat 2-dimensional quantum tori
as limits of elliptic curves appeared
in [CDS], but it was not made precise there.
Motivated by the preliminary version of present paper
Yu. Manin suggested in [M2] interesting ideas about quantum tori and abelian
varieties defined over arbitrary complete normed fields.}

\subsection{}

Relevance of derived categories (more technically,
triangulated $A_{\infty}$-categories
with finite-dimensional Hochschild cohomology) can be roughly
justified by the
following arguments:

1) Deformation theory of an $A_{\infty}$-category is controlled
by its Hochschild complex. Being defined properly (see [KoSo1]) such a category
always appears together with the formal moduli space of its deformations.

2) Reconstruction theorems, mainly due to A. Bondal, D. Orlov and
A. Polishchuk. Here are two examples.

\begin{thm} ([BO]). Let $X$ be a smooth irreducible variety with ample 
canonical or anticanonical sheaf, and $Y$ be a smooth
algebraic variety. If the bounded derived categories
of coherent sheaves $D^b(X)$ and $D^b(Y)$ are equivalent 
 as triangulated categories, then $X$ is isomorphic to $Y$.

\end{thm}

If $X$ is a Calabi-Yau manifold (for example, an elliptic curve)
then the theorem is not true. Nevertheless, often one can recover
from $D^b(X)$ some information about $X$. 

\begin{thm} ([O1]). For each abelian variety $X$ there are finitely many
non-isomorphic abelian varieties which have the
bounded derived category of coherent sheaves equivalent to $D^b(X)$.

\end{thm}

\subsection{}
Physicists discovered quantum tori in
various theories ( see [CDS], [SW]).
Morita equivalence of quantum tori was interpreted as a new duality
(see [RS], [S1]).
One hopes  that the
 concept of non-commutative
compactification will help to construct and
investigate mathematical models for dualities in 
quantum physics.

As an example  let us consider the Homological Mirror Conjecture (HMC) of
Kontsevich (see [Ko1]).
Homological Mirror Conjecture says that
the category $D^b(E_q)$ is equivalent to the bounded
 derived category $D^b(F(T^2_q))$
of the Fukaya category $F(T^2_q)$ (see for ex. [Ko1] for
the definition of the Fukaya category) of the symplectic torus
 $T^2_q=(T^2, 1/Re(\tau) dx_1\wedge dx_2)$. 
Let us imagine, that as $q$ approaches to a point
at the circle $|q|=1$, both derived categories degenerate into the same 
$\A$-category. Hence non-commutative
degenerations become manifestly equivalent. 
This might give an insight to HMC.
 We will
speculate about the category that can appear as such degeneration.
 Roughly speaking,
it is the (derived) category of bundles with connections which are flat
along a foliation. The foliations appear as degenerate
complex structures. Global sections of such
bundles are modules over the algebra of functions on a quantum torus. 
Thus one gets a
``non-commutative" description of the corresponding boundary stratum
of the compactified universal covering 
of the moduli space of complex structures.

The foliations which appear in the story are not arbitrary.
 They carry affine structures on leaves. This
makes the whole picture similar  to the topological mirror symmtery
of Strominger-Yau-Zaslow (see [SYZ]). In [SYZ] Calabi-Yau manifolds
are foliated by special Lagrangian tori, hence fibers carry affine structures.
On the other hand, the SYZ-picture is related to the ``large" complex
structure limit, while our approach seems to be of different nature.
Hopefully, they  both 
correspond to two different strata
of the boundary of the compactified moduli space
of $N=2$ superconformal field theories. The commutative
stratum  was discussed in [KoSo2] in the case of abelian varieties
(see also Section 5 below).

\subsection{}
Two-dimensional quantum tori can be interpreted in many different ways.
For example, they are Morita equivalent to
algebras of foliations. 
Quantum tori can also be thought of as quantizations of tori equipped with
 constant symplectic structures (they appear in the open string theory in this
way).
Hence they  fit into the
 framework of deformation theory. Deformation theory of a Poisson manifold
 $X$ gives rise to a formal family of categories 
 $\widehat{{\cal C}}_{\hbar}$
 of modules over a quantized algebra $C^{\infty}(X)_{\hbar}$ of smooth function
 on $X$, where $\hbar$ is the formal parameter. The problem is to
 construct a ``global" moduli space. The germ at ${\hbar}=0$ of this
 moduli space contains
 families of $A_{\infty}$-categories ${\cal C}_{\hbar},
 \hbar\in {\bf R}$ such that the formal completion of ${\cal C}_{\hbar}$
 at $\hbar= 0$ gives $\widehat{{\cal C}}_{\hbar}$. In general ${\bf R}$ can be
 replaced by some parameter space ${\cal M}_X$. If the global parameter space exists,
 one can speak about ``dualities" for the ``family of theories parametrized
 by ${\cal M}_X$". 
  In the case of quantum tori,
 the dualities can correspond 
 to Morita equivalence of quantum tori, or to equivalence of
some subcategories of the full categories of modules.
 The ``duality group" should be $SL(2,{\bf Z})$ in order
to agree with the interpretation of quantum tori as
 degenerate elliptic curves.
 
For general Poisson manifolds one can ask the following question.

 \begin{que} For a given Poisson manifold $X$, is there a 
 global parameter
 space ${\cal M}_X$? If it exists, how to describe points
 corresponding to the equivalent categories? What is the ``duality
 group"  (and why it is a group)?

 \end{que}
 
 We are going to discuss later in the paper
conjectures motivated by quantum groups. Interesting ideas 
in this direction can be found in [Fad].

\subsection{}

Here is the content of the paper. In Section 2 we discuss analogies
between quantum tori over ${\R}$ and abelian varieties over ${\C}$.
In Section 3 we explain why quantum tori can be thought of as 
points of the boundary of the compactified Teichm\"uller space
of an elliptic curve. Main idea here is to treat the Thurston boundary
(given geometrically in terms of foliations) as a non-commutative
space (i.e. as a category). Section 4 is devoted to mirror symmetry
for quantum tori. We suggest to ``compactify" the conventional
homological mirror symmetry, extending it to the boundary of the
universal covering of the moduli space of complex structures
(Teichm\"uller space in the case of elliptic curves).
Section 5 contains a comparison with the
``commutative" picture. This material is borrowed from [KoSo2].
In Section 6 we discuss possible relations to theta functions and
quasi-modular forms. We end this section with 
conjectures about dualities in quantum groups.

{\it Acknowledgements.} 
 I would like to thank A. Bondal,
A. Connes, P. Deligne, M. Douglas, L. Faddeev, V. Ginzburg,
 K. Gawedzki,  M. Gromov,
 D. Kazhdan,
Yu. Manin, W. Nahm,
D. Orlov, A. Polishchuk,
A. Rosenberg, A. Schwarz, 
 D. Zagier for useful discussions. 
I am grateful to Maxim Kontsevich, 
who shared with me many of his ideas 
and supplied with illuminating examples. Some of them
are used in the paper (especially in Sections 2.3, 3.1, 3.2, 5).
I thank to S. Barannikov for the comments on the paper.

I also thank IHES for hospitality 
and CMI for financial support. 
 
\section{Quantum tori and their representations}

\subsection{Generalities}

We start with the definition 
of an algebraic quantum torus (see for example [M1]).
Let  $K$ be a complete normed field . 

\begin{dfn}
Quantum torus $T(L,\alpha)$ of rank $d$ (or dimension $d$)
 is defined by a free ${\bf Z}$-module $L$
of finite rank $d$ (lattice) and a  bilinear form $\alpha:
L\times L\to K$ such that

$$ 
\alpha(m,n)\alpha(n,m)=1, \alpha(m+n,l)=\alpha(m,l)\alpha(n,l).
$$
More precisely, the coordinate ring $A(T(L,\alpha))$ of the quantum torus 
 is a $K$-algebra with the unit generated
by generators $e(n), n\in L$ subject to the relations

$$
 e(m)e(n)=\alpha(m,n)e(m+n).
$$

\end{dfn}

Algebra of analytic functions $A^{an}(T(L,\alpha))$ is a completion
of $A(T(L,\alpha))$ which consists of formal series
$\sum_{n\in L}a_ne(n)$ with the coefficients $a_n$ decreasing
as $|n|\to \infty$ faster than any power of $|n|$.

Assume that $K$ 
carries an involution $x\mapsto \overline {x}$, such that $|x\overline x|=|x|^2$.
Suppose that 
$\alpha(m,n)^{-1}=\overline{\alpha(m,n)}$ for any $m,n \in L$.
Then there is a natural 
involution of $A(T(L,\alpha))$ given by $e(m)^*=e(-m)$.
It makes $A_{\alpha}=A(T(L,\alpha))$ into  a $K$-algebra with involution.
Clearly in this case $\alpha$ takes value in the subgroup
$K_1=\{x\in K, |x|=1\}$. We will call such quantum tori
{\it unitary}. We consider their coordinate rings as objects in the
category of algebras with involutions. 

\begin{dfn} Let $T(L,\alpha)$ be a unitary quantum torus. The algebra
$B_{\alpha}:=C^{\infty}(T(L,\alpha))$ of smooth functions on $T(L,\alpha)$
consists of series $f=\sum_{n\in L}a_ne(n)$ where 
the sequence $a_n$ decreases faster than any power of $|n|$, as
$|n|\to \infty$.

\end{dfn}

This terminology is justified by the case $K={\C}$, where
one gets the algebra of smooth functions on a real torus.
One has the natural embedding of algebras with involutions: 
$A_{\alpha}\subset B_{\alpha}$. 

 We are going to consider unitary quantum tori over the field ${\C}$
 unless we say otherwise. The corresponding classical
torus $L_{\bf R}/L$ will be denoted by $T(L)$
(here $L_{\bf R}=L\otimes {\bf R}$). Let $\Phi:L_{\bf R}\to L_{\bf R}^{\ast}$
be a linear map such that $\Phi(x)(y)+\Phi(y)(x)=0$.
We set $\varphi(x,y)=\Phi(x)(y)$. Thus we get a skew-symmetric bilinear
form $\varphi:L_{\bf R}\times L_{\bf R}\to {\bf R}$.
Taking $\alpha(x,y)=exp(2\pi i \varphi(x,y))$ we obtain a
unitary quantum
torus, which will be denoted by $T(L,\varphi)$ or $T(L,\Phi)$.
The corresponding algebras of functions (algebraic and smooth)
will be denoted either by $A_{\varphi}, B_{\varphi}$ or by 
$A_{\Phi}, B_{\Phi}$.

Let $V_L=L_{\bf R}\oplus L_{\bf R}^{\ast}$. We equip this vector space
with an {\it odd} symplectic structure. It is given by the canonical
symmetric form $Q((x_1,l_1),(x_2,l_2))=l_1(x_2)+l_2(x_1)$. Then the
skew-symmetric map $\Phi$ defines a Lagrangian subspace 
$graph(\Phi)\subset V_L$. 

One can define the Grassmannian $Gr_0(V_L)$ of Lagrangian subspaces
in $V_L$. Then the map $\Phi\mapsto graph(\Phi)$ identifies quantum
tori with an open subset of $Gr_0(V_L)$. Let $O(V_L,Q)$ be the group
of linear automorphisms of $V_L$ preserving the form $Q$ (orthogonal group),
and $SO(V_L,Q)$ be the corresponding special linear group. Then
$O(V_L,Q)$ and $SO(V_L,Q)$ act
 transitively on $Gr_0(V_L)$. We will denote
by $SO(L,L^{\vee})$ the subgroup of $SO(V_L,Q)$ which preserves
the lattice $L\oplus L^{\vee}$, where $L^{\vee}=Hom(L,{\bf Z})$.

\subsection{Morita equivalence}

The following theorem was proved in [RS] in the framework
of $C^{\ast}$-algebras. 

\begin{thm} In the notation of the previous subsection, let $graph(\Phi)$ and 
$graph(\Phi^{\prime})$
are conjugate by an element of the group $SO(L,L^{\vee})$. Then
the algebras $B_{\Phi}$ and $B_{{\Phi}^{\prime}}$ are Morita equivalent.

\end{thm}

Choosing a basis in $L$ one can identify the group $SO(L,L^{\vee})$
with the group $SO(d,d, {\bf Z})$ of linear automorphisms of
the vector space ${\bf R}^{2d}$ preserving the form 
$\sum_{1\le i\le d}x_ix_{i+d}$ and the lattice ${\bf Z}^{2d}$.

Surprisingly, the same group appears in a different problem
concerning derived categories of coherent sheaves on complex abelian
varieties. We recall the following result.

Let $X$ and $Y$ be complex abelian varieties, $\widehat X$
and $\widehat Y$ are dual abelian varieties. We denote by $L_X$ and $L_Y$
the lattices of first homologies of $X$ and $Y$. 

\begin{thm} ([O1]).
The derived category $D^b(X)$ is equivalent
to $D^b(Y)$ iff there exists an isomorphism
 $X\times \widehat X\to Y\times \widehat Y$
which identifies $L_X\oplus L_X^{\vee}$ and $L_Y\oplus L_Y^{\vee}$
as odd symplectic lattices (both lattices are equipped with the canonical
symmetric forms $Q_X$ and $Q_Y$ as above). 
\end{thm}

\subsection{General algebraic scheme}

We refer the reader to [KoSo2] for the background on
$\A$-categories.
We are going to discuss here a general algebraic scheme which sheds some
light on the similarity between the derived category of coherent sheaves
on an abelian variety and the derived category of modules
over the algebra of functions on a quantum torus. 
We will impose some restrictions on the objects.
These restrictions can be relaxed. The reader should assume
that ``natural" conditions are imposed, so that ``everything works".

Let ${\cal C}$ be an $A_{\infty}$-category over
a field $k$ of characteristic zero, such that its Hochschild
cohomology $HH^i({\cal C})$ is finite-dimensional for all
$i\ge 0$. In what follows we will assume that $k={\bf C}$.

We assume that $HH^0({\cal C})$ is a $1$-dimensional
vector space. It can be thought of as a Lie algebra of the
group $Aut(Id_{{\cal C}})$ of automorphisms of the identity functor.
The total cohomology space 

$$\oplus_{i\ge 0}HH^i({\cal C})=\oplus_{i\ge 0}Ext^i(Id_{{\cal C}},Id_{{\cal C}})$$
carries a graded Lie algebra structure 
(with the Gerstenhaber bracket on the Hochschild cohomology).
Let ${\cal C}=D^b(X)$,
where $X$ is a smooth complex projective variety.
Then
$\oplus_{i\ge 0}HH^i({\cal C})=\oplus_{i,j\ge 0}H^i(X,\bigwedge ^{j}T_X)$
where $T_X=T^{1,0}_X$ is the holomorphic tangent bundle.

We assume that the Lie subalgebra $HH^1({\cal C})$ is abelian.
It admits the following interpretation. Let us
consider the group $Aut({\cal C})$ of automorphisms 
of the $A_{\infty}$-category ${\cal C}$. 
 It is the group of classes
of isomorphisms $[F]$ of equivalence functors $F:{\cal C}\to {\cal C}$.
 We assume that it carries a structure of a Lie group.
The Lie algebra $g_{\cal C}$
of the connected component of the unit of $Aut^0({\cal C}):=G_{\cal C}$ 
is isomorphic to
$HH^1({\cal C})$. Under our assumptions, the group $G_{\cal C}$ is a 
finite-dimensional commutative Lie group over ${\bf C}$.

There exists a bundle $P$ over 
$G_{\cal C}\times G_{\cal C}$ with the fiber which 
is a ${\bf C}^{\ast}$-torsor. Let us describe it in detail.
 Let $F$ and $H$ be two functors, such that their isomorphism 
classes $[F],[H]$ belong to the group
$G_{\cal C}$. 
 According to our assumption on $HH^0({\cal C})$, the set
of isomorphisms $Iso(F\circ H,H\circ F)$ is a ${\bf C}^{\ast}$-torsor.
Suppose that $F_1$ and $H_1$ are another representatives of
the classes $[F]$ and $[H]$ respectively. Then $F\sim F_1$ and
$H\sim H_1$. This gives rise to an isomorphism of torsors
$Iso(F\circ H,H\circ F)\simeq Iso(F_1\circ H_1,H_1\circ F_1)$.

The following lemma is easy to prove.

\begin{lmm} This isomorphism of torsors depends on the
equivalence classes $[F]$ and $[H]$ only.
In other words, for different liftings of $[F]$ and $[H]$, we
can canonically identify the torsors.

\end{lmm}

Thus  we obtain a ${\bf C}^{\ast}$-torsor $P=P_{\cal C}$ 
over  $G_{\cal C}\times G_{\cal C}$
with the fibers 
$P_{[F],[H]}=Iso(F\circ H,H\circ F)$.

The abelian group $L_{\cal C}:=H_1(G_{\cal C},{\bf Z})$
carries a bilinear form 
$(,):L_{\cal C}\times L_{\cal C}\to {\bf Z}$, such that
$(x,y)=c_1(P)(x\boxtimes y)$, where $c_1(P)$ is the first Chern class
of $P$. We will keep the same notation for the
${\bf C}$-linear extension of the bilinear form to $L_{\cal C}\otimes {\bf C}$.

Since $G_{\cal C}$ is a Lie group, its fundamental group
is commutative, and hence it is isomorphic to $H_1(G_{\cal C},{\bf Z})$.
Let $\gamma$ be a loop based at the unit $e$ of $G_{\cal C}$.
We define a linear map $p:L_{\cal C}\otimes {\bf C} \to g_{\cal C}$ 
by the formula $p(\gamma):=\dot {\gamma}(e)$ 

Then we have the 
following exact sequence of linear maps:

$$0\to \Lambda_{\cal C}\to L_{\cal C}\otimes {\bf C} \to g_{\cal C},$$
where $\Lambda_{\cal C}$ is defined as the kernel of $p$.

\begin{conj} The symmetric bilinear form $(x,y)$ is
non-degenerate, and the subspace $\Lambda_{\cal C}$ is maximal isotropic 
(i.e. Lagrangian) with respect to it.

\end{conj}

It is not clear how to prove this conjecture for general
$A_{\infty}$-categories.
In two examples considered below (abelian varieties and
quantum tori) the proofs are straightforward.

Assuming the conjecture, to the $A_{\infty}$-category ${\cal C}$ 
we have assigned canonically
a Lagrangian subspace $\Lambda_{\cal C}$ in the (odd) symplectic
vector space $L_{\cal C}\otimes {\bf C}$.
It is interesting to notice that the linear algebra data seem
to contain certain information about ``discrete" symmetries of
the $A_{\infty}$-category. Apriori this was not obvious.
\begin{rmk}
In the case of quantum tori we will need to consider
non-Hausdorff Lie groups. Then the abelian Lie algebra $L_{\cal C}$
 should be defined
not as the first homology, but as the kernel of the exponential
map $exp: g_{\cal C}\to G_{\cal C}$. Notice that $g_{\cal C}$ is
defined canonically: it is the Lie algebra of infinitesimal symmetries of our
$A_{\infty}$-category. The group $G_{\cal C}$ is not defined
canonically. At the level of local
Lie groups one can think about $G_{\cal C}$ as about quotient
of the local Lie group corresponding to $g_{\cal C}$ by the discrete
subgroup consisting of such $\gamma\in g_{\cal C}$ that $exp(\gamma)$
gives rise to an automorphism of the category ${\cal C}$ isomorphic
to the identity functor. 
\end{rmk}

\subsection{Example: abelian varieties and quantum tori} 

Let ${\cal C}=D^b(X)$ where $X$ is a complex abelian variety
of dimension $n$. Then $G_{\cal C}\simeq {\bf C}^{2n}/{\bf Z}^{4n}$.
The bundle $P$ is basically the tensor square of the Poincare line
bundle on $X\times \widehat{X}$ (this follows from the results of
Polishchuk and Orlov, see [O1]).
The odd symplectic form is the
canonical symmetric form on
 ${\bf C}^{4n}={\bf C}^{2n}\oplus ({\bf C}^{ 2n})^{\ast}$.
It is easy to see that 
$HH^{\ast}({\cal C})\simeq \bigwedge^{\ast}({\bf C}^{2n})$. The latter
can be also interpreted as the algebra of functions on the odd Lagrangian
subspace in ${\bf C}^{4n}$. Thus the formal deformation theory of ${\cal C}$ 
as an 
$A_{\infty}$-category is the same as the deformation of the corresponding
odd Lagrangian subspace as a (non-linear) Lagrangian submanifold 
$\Lambda_{\cal C}\subset {\bf C}^{4n}$.

For a quantum torus $T(L,\Phi)$ of dimension $n$ one takes as ${\cal C}$
the derived category of modules of finite rank over the algebra 
$B_{\Phi}(L)$. It is expected that automorphisms of the category come
from automorphisms of the algebra.
Then $G_{\cal C}\simeq T^n/{\bf Z}^n$ is a non-Hausdorff
Lie group (it is the group of automorphisms of the algebra $B_{\Phi}(L)$
modulo inner automorphisms). Following the general scheme outlined above,
one obtains a Lagrangian subspace 
(it coincides with $graph(\Phi)$) in the odd symplectic
vector space ${\bf R}^{2n}={\bf R}^n\oplus ({\bf R}^n)^{\ast}$.
For the categories discussed in this subsection
 the cohomology $HH^0$ is $1$-dimensional, and $HH^1$ is
a commutative Lie algebra. Moreover, the Conjecture $1$ holds in both cases.

\section{Complex structures and foliations}

\subsection{Foliations and degenerate complex structures}

Let $X$ be an even-dimensional smooth real manifold, which carries
a complex structure. The latter is given by an integrable
 subbundle $T^{0,1}_X\subset T_X\otimes {\bf C}$ of anti-holomorphic directions,
where $T_X$ is the real tangent bundle of $X$. 
Suppose that we have a family of
complex structures degenerating into a real foliation
 of the rank equal to $dim_{\C}\,X$. We will call such degenerations
 {\it maximal}.
In fact one should consider 
subsheaves of the tangent sheaf
$T_X$ because the foliation can be singular.
Degeneration gives rise not just to a foliation $F$ of $X$, but also to
an isomorphism $j$ of vector bundles $F$ and $T_X/F$.
The isomorphism satisfies some integrability conditions, which
give rise to an affine structure on the fibers of $F$. 
To be more precise, the space $Hom(F,T_X/F)$ can be interpreted as
a tangent space $T_{F}({\cal M})$ to $F$ in the
``moduli space'' ${\cal M}$ of all foliations
on $X$. When $T^{0,1}_X$ and the holomorphic subbundle
 $T^{1,0}_X$ get close to each other, one
 obtains  a tangent vector in the space 
 $T_{F}({\cal M})$.
 Similarly, for any fiber $F_y,y\in X$ of the foliation $F$, the vector
space $(T_X/F)_y$ can be identified with the tangent
 space to the leaf $F_y$ in the ``moduli space''
space of all leaves of the foliation $F$
 (it does not depend on a point of the leaf). It is easy
 to check that the integrability condition implies the
 following result.

\begin{prp} The isomorphism $j:F\to T_X/F$ 
gives rise to an isomorphism
of the tangent space to a leaf $F_y$ in the space of leaves,
 with the space of commuting vector fields
on the leaf. In particular it yields an affine structure on the leaf $F_y$.

\end{prp}
Informally this result can be explained in the following way.
Tangent space
to a ``degenerating" complex manifold
is invariant with respect to the rotation
by 90 degrees (multiplication by $i=\sqrt{-1}$).
 Let us move a point $x$ along
some leaf $F_y$ of the limiting foliation $F$.
Then the rotation by 90 degrees of vectors
from $T_xX/F_x$ transforms them into $F_x$. On the other hand,
the bundle $T_X/F$ carries a canonical flat connection (Bott connection).
We can identify the spaces $T_{X,x}/F_x$ 
for close points $x$ using the connection. 
Hence a choice of commuting vectors at the space $T_{X,y}/F_y$
gives rise to commuting vector fields along the leaf $F_y$.

One can show that the space of pairs $(F,j)$ 
being factorized by
the natutal action of the multilicative group ${\bf R}^{\ast}$
(dilations of $j$) is generically a subspace of the real codimension $1$
in the space of integrable subbundles in $T_{\bf C}X$.
The latter can be considered as a subvariety of the total
space of the bundle of  Grassmannians  $\Gamma(X,Gr(T_{\bf C}X))$. 
In this way one gets a compactification of the universal covering
of the moduli space of complex structures on $X$. In the
case of curves the compactification adds the
Thurston boundary to the Teichm\"uller space.

One can try to ``compactify" the
category of coherent sheaves on a complex manifold.
The category of 
sheaves equipped with flat connections along the
 foliation, which is the maximally degenerate complex structure,  serves as a ``point'' of the
``non-commutative'' stratum of the boundary.
In the case of elliptic curves the compactification is compatible
with the natural $SL(2,{\bf Z})$-action. Indeed, in both cases
(complex structure given by $\tau\in {\cal  H}$ and affine foliation on $T^2$
given by $dt=\varphi dx$) if the values of parameters are
$SL(2,{\bf Z})$-conjugate, then there is an automorphism of $T^2$
(as a smooth manifold) which identifies the corresponding subbundles
(in $T_{\bf C}T^2$ and $T_{\bf R}T^2$ respectively). 

More generally, one considers the space of complex polarizations on
a given compact complex manifold $X$. To every polarization
$\tau$ one assigns the corresponding bounded derived category
of coherent sheaves. The space of polarizations admits a natural 
compactification by real foliations, as we discussed above.
The question is: what is the category which should be assigned
to a foliation $F$ which is a point of the boundary? In the case of maximal degeneration we suggest
to take the (derived) category of $F$-local systems (see next subsection).

\begin{rmk} Compactification by pairs $(F,j)$
 modulo the action of ${\R}^{\ast}$ is not compatible 
with the action of the group of diffeomorphisms $Diff(X)$. In particular,
it does not descends to a compactification of the moduli space
of complex structures on $X$. 
In 1-dimensional case the action of the mapping class group
extends to the Thurston boundary of the Teichm\"uller space.

\end{rmk}

\subsection{Foliations and ``small" modules over quantum tori}

Let $X$ be a smooth manifold of dimension $2n$, and $F$ be a foliation
of $X$ of rank $n$. Let $W$ be a sheaf of $C^{\infty}_X$-modules on $X$,
where $C^{\infty}_X$ is the sheaf of smooth functions on $X$.

\begin{dfn} We say that $W$ carries an $F$-connection, if we are given
a morphism of sheaves $\nabla: W\to W\otimes F^{\ast}$ such that
$\nabla_v(fs)=f\nabla(s)+v(f)s$ for any germs $f\in C^{\infty}_X$,
$s\in W$, and $v\in F$ (here we identify the subbundle $F$ with the
sheaf of sections), $F^{\ast}=Hom(F,C^{\infty}_X)$.

\end{dfn}

The category $Sh(X,F)$ of sheaves of finite rank, which carry an
 $F$-connection
form a tensor category. The curvature of an $F$-connection
is defined in the usual way. Sheaves which carry $F$-connections with
zero curvature are called $F$-flat. If an $F$-flat sheaf is locally
free (i.e. it corresponds to a smooth vector bundle on $X$), it
is called an $F$-local system. 
We denote by $D^b(X,F)$ the bounded derived category of
the category $Loc(X,F)$ of $F$-local systems
on $X$. If $X$ carries
a symplectic form $\omega$, and $F$ is a Lagrangian foliation then
 we will sometimes add $\omega$ to the notation. 
Foliations described in the Proposition 1 will be of the main interest for us.
Here is an example.
Let $F$ be an affine foliation of rank $n$ on the standard torus 
$T^{2n}={\bf R}^{2n}/{\bf Z}^{2n}$. This means that $F$ is defined
as $(V\oplus {\bf R}^{2n})/{\bf Z}^{2n}$
for some $n$-dimensional vector subspace
$V\subset {\bf R}^{2n}$. Let us choose the subspace
$S$ such that $S\oplus V={\bf R}^{2n}$,
and $S$ defines a closed $n$-dimensional submanifold $Y$ in $T^{2n}$. 
There is a pull-back functor from the category $Loc(T^{2n},F)$ to the category
of vector bundles on $Y$. 
On the other hand, if $(\Lambda,\nabla)$ is an $F$-local system
on $X$, then the holonomy of the connection $\nabla$ defines an action of
the group ${\bf Z}^n$ on the restriction $\Lambda_{|Y}$. 
Since ${\bf Z}^n$ acts on $Y$,
we obtain a structure of an ${\bf Z}^n\ltimes C^{\infty}(Y)$-module
on the space of section $\Gamma(Y,\Lambda)$. We can complete the 
cross-product
algebra thus getting the algebra of smooth functions on a 
quantum torus acting on $\Gamma(Y,\Lambda)$.

\begin{dfn} a) Let $B_{\Phi}(Y)$ be the algebra of smooth
 functions on the quantum
torus described above. We call a $B_{\Phi}(Y)$-module {\it small}
if it is projective as a module over the commutative subalgebra
$C^{\infty}(Y)$.

b) More generally, let $B_{\Phi}$ be an algebra of smooth functions 
on a quantum torus $T(L,\Phi)$ such that $rk(L)=2n$, and $\Phi$ defines
a symplectic structure on $L_{\bf R}=L\otimes {\bf R}$.
 Let $L_0\subset L$ 
be a Lagrangian sublattice (i.e. $L_0\otimes {\bf R}$ is a Lagrangian
subspace in   $L_{\bf R}$). We say that a $B_{\Phi}$-module $M$
is small with respect to $L_0$, if it is projective with respect
to the maximal commutative subalgebra of  $B_{\Phi}$ spanned
by $e(\lambda),\lambda\in L_0$.

\end{dfn}

Clearly small modules with respect to a given
Lagrangian sublattice $L_0$ form a category $B_{\Phi}(L_0)-mod$.
Let us consider an example of $2$-dimensional tori. Then $rk(L)=2$,
$rk(L_0)=1$. We will identify $L$ with ${\bf Z}^2$ and $L_0$ with
${\bf Z}\oplus 0\subset {\bf Z}^2$.
 Let $\varphi\in {\bf R}\setminus {\bf Q}$, and
$F$ be the affine foliation $dt=\varphi dx$ in the standard coordinates
in ${\bf R}^2=L\otimes {\bf R}$. 

\begin{prp} The category $B_{\varphi}(L_0)-mod$ is equivalent to the
category of $F$-local systems on the torus $T(L)=L_{\bf R}/L$.

\end{prp}

{\it Proof}. We have already constructed a functor from $F$-local systems
to the category of small modules. Namely, every $F$-local system $V$
being restricted
to the equator of $T^2$, gives rise to a projective module over the algebra
$C^{\infty}(T^1)$. Since $\varphi$ is irrational, the foliation defines
an action of the group ${\bf Z}$ on $T^1=(x\, mod\,{\bf Z},0)$. The holonomy
of the flat connection defines the structure of a small module
on $\Gamma(T^1,V_{|T^1})$.

An inverse functor is constructed such  as follows. Let $M$ be a small module
over $B_{\varphi}$. We denote the standard generators of $B_{\varphi}$
by $e_1$ and $e_2$. Then $M$ is a projective module over the subalgebra
$B_{\varphi}(e_2)$ generated by $e_2$. Thus we have a vector bundle
$\widehat{M}$ over $T^1=\{(x\, mod\, {\bf Z}, 0)\}\subset T^2$, such that
 $\Gamma(T^1,\widehat{M})=M$. 
 
Let $V\subset  \Gamma(T^1,\widehat{M})\otimes C^{\infty}({\bf R})$ 
consists of 
elements $f$ such that
$f(x+1,t)=f(x,t)$ and $f(x,t+1)=e_1(f(x-\varphi,t))$ ( here we write formulas
in coordinates $(x,t)\in {\bf R}^2$ rather than in coordinates
$(x\, mod\, {\bf Z}, t)\in T^1\times {\bf R}$). 
It is easy to check that they are global sections of the vector bundle
$V\to T^2$ such that for its pull-back to
${\bf R}^2$ we have: the fiber $V_{(x,t)}$ is naturally isomorphic to
 $\widehat{M}_x$. We define an
$F$-connection $\nabla_F$ on $V$ by identifying infinitesimally 
closed fibers:
$g_{\varepsilon}:f(x,t)\mapsto f(x+\varepsilon \varphi, t+\varepsilon)$.
Then $g_1:f(x,t)\mapsto f(x+\varphi,t+1)=e_1(f(x-\varphi+\varphi,t))=
e_1(f(x,t))$. The action of the 
holonomy of $\nabla_F$ (shift on the period
$t=1$) is equivalent to the action of the generator $e_1$ on the
module $M$. Thus the action of ${\bf Z}\ltimes C^{\infty}(T^1)$ on
$\{f(x,0)\}=M$ given by the $F$-local system is the same as
the structure of $B_{\varphi}$-module on $M$. 
The Proposition is proved.$\blacksquare$

\begin{rmk} 1) In order to identify the (derived) categories of coherent sheaves on
the elliptic curves
$E_{\tau_1}$ and $E_{\tau_2}$ we use a bimodule, which is a 
sheaf of regular functions on the graph of an automorphism $f:T^2\to T^2$ identifying
 the complex structures.
If $p_i:E_{\tau_1}\times E_{\tau_2}\to E_{\tau_i}, i=1,2$ 
are natural projections, then the equivalence functor is given by
$M\mapsto p_{2\ast }(p_1^{\ast}M\otimes {\cal O}_{graph(f)})$.
In the case of tori with affine foliations we basically use the
same description. Notice that it does not give a Morita equivalence
of the corresponding quantum tori. It gives an equivalence
of the categories of small modules.

2) One should notice that vector bundles with $F$-connections
form a tensor category, while projective (or all finite)
modules over quantum tori do not. 

3) Small modules are similar to holonomic $D$-modules.
Their algebraic version was studied from this point of view in [Sab].

\end{rmk}

Suppose  we have a free abelian Lie group $G$
together with a dense embedding of $G$ into the group $Aut(T^n)$ of
the affine automorphisms of the torus $T^n$. Then to have a vector
bundle $V$ over $T^n$ together with a lifting of the action
of $G$ to $F$, is the same as to have a module over the quantum
torus defined as a completed cross-product of the group
algebra of $G$ and $C^{\infty}(T^n)$. We used this simple
observation in the case when the action of $G={\bf Z}$ was
induced by an affine foliation of $T^2$. In general we have
a functor {\it \{Affine foliations\}} $\to$ {\it \{Small modules\}}.
Notice that for $n>1$ there is no inverse to this functor.
In other words we cannot 
recover a foliation from a small module over the algebra
$B_{\varphi}$.
Indeed, we have less affine foliations than quantum tori.
Quantum torus is given by a skew-symmetric form. Hence the dimension of the
moduli space of quantum tori of the rank $2n$ 
is equal to $n(2n-1)$. At the same time the moduli space of affine foliations
has the dimension $n^2$ (every such a foliation is given by the graph
of a linear morphism ${\bf R}^n\to {\bf R}^n$).
 These two numbers coinside only if $n=1$.
 
 \subsection {Coherent sheaves on elliptic curves and quantum tori}
 
 Here we recall a result from [BG] and use it to provide another link
 between the derived category of coherent sheaves on elliptic curves
and the derived
 category of certain modules over quantum tori. Assume that we have
 fixed a non-zero complex number $q$ such that $|q|<1$.
Let us consider a complex algebra $\bar{A}_q$ which is generated by
the field of  Laurent formal power series ${\bf C}((z))$ and an
invertible element $\xi$ such that $\xi\,f(z)=f(qz)\,\xi$ for any
$f\in {\bf C}((z))$. One defines an abelian category ${\cal M}_q$
such as follows. Objects of ${\cal M}_q$ are $\bar{A}_q$-modules $M$, which are
${\bf C}((z))$-modules of finite rank. In addition, it is required that
there exists a free ${\bf C}[[z]]$-submodule $M_0\subset M$
of maximal rank, such that
$M_0$ is invariant with respect to the subalgebra 
${\bf C}[\xi,\xi^{-1}]\subset \bar{A}_q$. Clearly ${\cal M}_q$ is a
${\bf C}$-linear rigid tensor category.
It is proved in [BG] that ${\cal M}_q$ is equivalent to the
tensor category $Vect^{ss}_0(E_q)$ of  degree zero semistable
holomorphic vector bundles on the elliptic curve 
$E_q={\bf C}^{\ast}/q^{\bf Z}$. Let us recall the idea of the proof. Using the
Fourier-Mukai transform, applied to the vector bundles in question,
one gets the category of sheaves with finite support. The latter
category admits a description in terms of linear algebra data
(a vector space equippped with an endomorphism). The same data
describe objects of ${\cal M}_q$.

If we forget about tensor structures, then
 the category  $Vect^{ss}_0(E_q)$ generates 
a subcategory$D^b_{fin}(E_q)$ of $D^b(E_q)$. 
Here we understand
the word ``generate" in the following sense:
one allows to take extensions, direct summands of objects
and all shifts of an object. We will use the notation $D^b({\cal M}_q)$
for the derived category generated by ${\cal M}_q$.
Then the previous discussion implies the following result.

\begin{prp} The categories $D^b({\cal M}_q)$ and $D_{fin}^b(E_q)$
are equivalent.

\end{prp}

One can use this theorem as a {\it definition} of $D_{fin}^b(E_q)$
in the case when $|q|=1$.
\begin{que} Let $|q|=1$. Is
it true that $D^b({\cal M}_q)$ is equivalent to a subcategory of the derived
category of small modules over the quantum torus $B_q$?

\end{que}

\begin{que} Let $|q|<1$.
How to describe the category $D^b(E_q)$ in terms of modules
over $B_q$ (or some related algebra)?

\end{que}

In order to answer the last
question, one can consider the algebra generated by
holomorphic functions on ${\C}^{\ast}$
and shifts $z\mapsto qz$. The question is whether it is possible
to replace the algebra of holomorphic functions on ${\C}^{\ast}$
by something ``more algebraic". Then one would have a uniform
description of the derived category of coherent sheaves on
an elliptic curve and its ``non-commutative" degeneration.

Hopefully, the above  considerations can be generalized
to the case when ${\C}((z))$ is replaced by a complete normed field
(for example, to the $p$-adic case). Then it can be used as a definition
of the derived category of coherent sheaves on quantum tori defined
over such fields  (cf. [M2]).

\section{Mirror symmetry and deformation quantization}

\subsection{Reminder on Homological Mirror Conjecture}

Homological Mirror Conjecture (HMC) was formulated by Kontsevich in 1993.
We are not going to recall all the details here (see [Ko1], [KoSo2]). 
HMC is a claim about equivalence of two triangulated
$A_{\infty}$-categories $D^b_{\infty}(X)$ and $F(X^{\vee})$
for given  mirror dual complex
Calabi-Yau manifolds $X$ and $X^{\vee}$. 
The category $F(X^{\vee})$ is the Fukaya category of $X^{\vee}$.
It can be defined for any symplectic manifold $(M,\omega)$.
Objects of $F(M)$ are pairs $(N,L)$ where $N$ is a Lagrangian submanifold
of $M$ and $L$ is a unitary local system on $N$. For
two objects $(N_1,L_1), (N_2,L_2)$ such that $N_1$
and $N_2$ intersect transversally, one defines the space
of morphisms as 

$$Hom((N_1,L_1), (N_2,L_2))=\oplus_{x\in N_1\cap N_2}Hom(L_{1x},L_{2x})$$
where $L_{ix}, i=1,2$ are fibers of the local systems.

The structure of $A_{\infty}$-category on $F(M)$ is given in terms of the following data:

1) A structure of complex on all spaces
$Hom((N_1,L_1), (N_2,L_2))$.
In particular they are ${\bf Z}$-graded complex vector spaces
with the grading defined by means of the Maslov index.

2) Higher compositions, which are morphisms of complexes

$$m_k:\otimes_{0\le i\le k}Hom(X_i,X_{i+1})\to Hom(X_0,X_{k+1})[2-k],$$

for given objects $X_i=(N_i,L_i)\in F(M)$. The composition $m_k$
is defined in terms of the moduli space
of holomorphic maps of $(k+1)$-gons $C_{k+1}$
to $M$, such that the $i$th side of $C_{k+1}$ belongs to $N_i$.
Thus, $m_1$ is basically the differential in the Floer complex
associated with the pair of Lagrangian submanifolds $N_1$ and $N_2$.
Formulas for $m_k,k\ge 2$ involve also the monodromies of flat connections
along the sides of polygons as well as areas of the polygons computed
with respect to the symplectic form on $M$.
There are compatibility conditions for the morphisms $m_k$.

The category $D_{\infty}^b(X)$ is an $A_{\infty}$-version of the derived
category of coherent sheaves on $X$. Its objects are 
bounded complexes of holomorphic vector bundles.
 Essentially the same $A_{\infty}$-category 
is given by the dg-category of dg-modules over the dg-algebra of Dolbeault
forms $\Omega^{0,\ast}(X)$.

\subsection{The case of elliptic curves}

In the case of elliptic curves the HMC was proved in [PZ] (see also [AP]).
One starts with a symplectic $2$-dimensional torus
$(T^2,\omega)$, where $\omega$ is a constant symplectic form.
In fact one should also fix a real number $B$ (more precisely,
the cohomology class 
$[Bdx\wedge dy]\in H^2(T^2,{\bf R})/H^2(T^2,{\bf Z}))=
{\bf R}/{\bf Z}$). Then one defines the Fukaya category  $F(T^2)$
using the complexified symplectic form $Bdx\wedge dy+i\omega$.
All the definitions are standard, only the areas of polygons
become complex numbers (see details in [PZ]). 
The symplectic form $\omega$ determines a unique flat metric
$A(dx^2+dy^2)$ on $T^2$ such that $A=\int_{T^2}\omega$.
It gives a unique complex structure on $T^2$ (because the size
of the chamber of the lattice ${\bf Z}^2$ is fixed).
In this way we obtain an elliptic curve (i.e. $1$-dimensional
Calabi-Yau manifold).
Let $X=E_{\tau}={\bf C}^{\ast}/q^{\bf Z}$ be this curve,
$q=exp(2\pi i \tau), Im(\tau)>0$.
Then the dual Calabi-Yau manifold is the elliptic curve
$E_{\rho}={\bf C}^{\ast}/e^{2\pi i\rho {\bf Z}},\rho=B+iA$.

Mirror symmetry functor can be described explicitly. On the
symplectic side of HMC one has closed $1$-dimensional submanifolds
in $T^2$ which carry unitary (quasi-unitary for the version
of HMC considered in [PZ]) local systems. On the complex
side of HMC one has complexes of holomorphic vector bundles
(they generate the derived category of coherent sheaves).
The dictionary between symplectic and complex sides
translates standard $(m,n)$ geodesics  equipped with trivial
$1$-dimensional local systems into holomorphic vector bundles
of rank $n$ with the first Chern class $m$.
For general abelian varieties a version of HMC was proved
in [KoSo2] using methods of Morse theory and non-arhimedean analysis.

\subsection{``Compactified" Homological Mirror Symmetry}

Homological mirror conjecture implies that formal moduli spaces
of deformations of  $D^b_{\infty}(X)$ and $F(X^{\vee})$ for dual
Calabi-Yau manifolds are isomorphic. 
The tangent space to the moduli space of deformations of an
$A_{\infty}$-category ${\cal C}$ is $\oplus_{k\ge 0}Ext^k(Id,Id)$. Here
one takes $Ext$ groups in the category of functors ${\cal C}\to {\cal C}$,
and $Id$ is the identity functor. In other words, the tangent space
is given by the Hochschild cohomology of the category. The Yoneda
product of $Ext$ groups gives rise to the product on
the tangent space.

In the case of $D^b_{\infty}(X)$ the tangent space is
$\oplus_{p,q\ge 0}H^p(X,\bigwedge^qT_X)$. In the case of
$F(X^{\vee})$ it is $H^{\ast}(X^{\vee})$. The product in
the latter case is the small quantum  cohomology product defined in terms 
 Gromov-Witten invariants (see [Ko1] for details).
The formal moduli spaces are bases of semi-infinite variations
of Hodge structures (see [B2]).
The mirror symmetry functor should identify the
corresponding Hodge filtrations
for $D^b_{\infty}(X)$ and $F(X^{\vee})$
(or, better, semi-infinite variations of Hodge structures,
as suggested by Barannikov, see [B3]).

  We mention here that
non-commutative analog of the variations of Hodge structures
with applications to mirror symmetry was 
introduced in [B2-B3].

In the case of elliptic curves (more generally, abelian varieties) there are no
 Gromov-Witten invariants, but the rest of the
picture is present. Then one can ask the following question.

\begin{que} What happens with both sides of HMC as $Im(\tau)\to 0$?

\end{que}

We have argued that the derived category of coherent sheaves on 
the elliptic curve $E_{\tau}$ degenerates into
the derived category of $F$-local systems,
where $F$ is a foliation on $T^2$ (degenerate complex structure).
Equivalently, it is the derived category of the category of small
modules over the algebra of functions on the corresponding quantum torus.
Another way to treat this degeneration is to consider the derived 
category of ``coherent sheaves on the non-commutative upper-half plane".
This means that we start with  an algebra $B$ which
is the algebra over ${\bf C}[q,q^{-1}]$ consisting
of series $f=\sum_{m,n}a_{m,n}x^my^n$ where $xy=qyx$, the coefficients
$a_{m,n}$ are rapidly decreasing, $m, n \in {\bf Z}$, and
$n$ runs through a finite set. If $|q|\ne 1$ modules 
of finite rank over this algebra
are the same as ${\bf Z}$-equivariant sheaves of finite rank 
on ${\bf C}^{\ast}$, i.e. the same as coherent sheaves on 
$E_q={\bf C}^{\ast}/q^{\bf Z}, q=e^{2\pi i \tau}$. If $|q|=1$
 then every such module
admits a structure of module over the algebra $B_{\tau}$ 
(i.e. we allow
infinite sums in $n$). Then the derived category of coherent sheaves
over $B$ can be treated as a family of derived categories, which
has as fibers  elliptic curves and quantum tori for different
values of $q$. Now the $SL(2,{\bf Z})$-invariance
of categories is manifest. In fact one has two copies of $SL(2,{\bf Z})$
acting on the family of categories. One copy acts by fractional transformations
of $\tau$. It gives rise to an isomorphism of elliptic curves if 
$Im(\tau)>0$, and  equivalence of the categories
of small modules if $Im(\tau)=0$. The equivalence in this case is given
by a bimodule $H_{\tau, g(\tau)}$,  where $g\in SL(2,{\bf Z})$.
One can check that $H_{\tau, g_1g_2(\tau)}\simeq 
H_{g_2(\tau), g_1g_2(\tau)}\otimes H_{\tau, g_2(\tau)}$,
and $H_{\tau, \tau}$ corresponds to the identity functor.

Another copy of $SL(2,{\bf Z})$ acts by automorphisms of the algebra
$B_q$ for a
fixed $q$ which is not a root of 1.
 Namely, $x\mapsto x^ay^b, y\mapsto x^cy^d$ is an algebra
automorphism iff $ad-bc=1$. Taking $a=d=0, b=-c=1$ one gets
an analog of the Fourier-Mukai transform in the case of quantum tori.

We treat this picture as a non-commutative compactification of the universal
covering of the moduli space of elliptic curves.

It is natural to expect that there is similar compactification
of the ``symplectic" side of HMC. In the next section we are going to present
some arguments in favor of the idea, that the degenerate Fukaya category
of the symplectic torus $(T^2,(1/\tau)\omega)$ is again the category of small
modules over the algebra $B_{\varphi}$, where $\varphi=Im(\tau)$, and
$B_{\varphi}$ is interpreted  as a quantized algebra
of functions on the torus.

The ``compactified" HMC should be a statement about an equivalence
of families of $A_{\infty}$-categories, both parametrized by
$\tau\in \overline {\cal H}=\{z\in {\bf C}|Im(z)\ge 0\}$. For $\tau \in {\cal  H}$
one has $A_{\infty}$-categories from the conventional HMC, and for 
$\tau \in {\bf R}$ one has two equivalent categories of modules for
the same quantum torus, but described in two different ways: 

a) as the category of $F$-local systems on $T^2$;

b) as the category of modules over a quantized symplectic torus, which
correspond to Lagrangian submanifolds.

Although two descriptions a) and b) give rise to equivalent categories,
the mirror symmetry functor produces a non-trivial identification
of the Hodge filtrations on the cohomology of $X=T^2$.
Hodge filtrations
are filtrations on the periodic cyclic homology of the $\A$-categories
under considerations.
Periodic cyclic homology of
either  $A_{\infty}$-category is
isomorphic to the total cohomology of the underlying space.
Mirror symmetry functor should 
interchange the Hodge filtrations on periodic cyclic homology. 

In the case of the 2-dimensional quantum torus $T^2({\bf Z}^2,\varphi)$
 one has two filtrations
on the periodic cyclic homology (the latter is isomorphic to the cohomology
of the usual torus $T^2$). First filtration arises from the interpretation
of the quantum torus as a quantized symplectic manifold. It has
one non-trivial term corresponding to the line spanned by the vector
$(1,[\omega])$, where $1\in H^0(T^2)$ is the unit in cohomology, and $[\omega]$
is the cohomology class of the symplectic form
 $\omega={1\over \varphi}dx\wedge dt$.
Second filtration arises from the interpretation of the quantum torus
in terms of the foliation. The only non-trivial term of this filtration 
corresponds to the straight
line spanned by the class of the foliation 1-form $dx-\varphi dt$.
Now we have filtrations in even and odd cohomology respectively.
The mirror symmetry functor interchanges odd and even cohomology,
interchanging also the Hodge filtrations.
This picture is a ``limiting" one for the homological
mirror symmetry in the case of elliptic curves. Therefore
the mirror symmetry functor in the case of quantum tori identifies the (derived) category
of modules over the algebra $B_{\varphi}$ of smooth functions
on $T(L,\varphi)$ with itself. One can say that the quantum torus is mirror dual
to itself. The mirror functor acts non-trivially interchanging odd and even cohomologies,
and identifying the Hodge filtrations described above
(cf. with the case of complex abelian varieties considered in  [GLO]).

\subsection{Modules over quantized algebras and Fukaya category}

In order to define the Fukaya category one needs a symplectic
manifold. There is a simpler abelian category depending
on symplectic structure. Let $(M,\omega)$ be a symplectic manifold.
Then it admits a deformation quantization (i.e. formal family
of associative products with the prescribed first jet).
 The quantized algebra
$C^{\infty}(M)_{\hbar}$ is not defined canonically. On the other hand,
the abelian category ${\cal C}_M$
 of $C^{\infty}(M)_{\hbar}$-modules is defined canonically.

Let $\Lambda\subset M$ be a Lagrangian submanifold. Identifying
a neighborhood of $\Lambda$ with a neighborhood of some manifold
$X$ in $T^{\ast}X$, one can construct a $C^{\infty}(M)_{\hbar}$-module
$W_{\Lambda}$, 
which is canonically defined as an object of ${\cal C}_M$. 
For the modules $W_{\Lambda}$ one expects
the following theorem (cf. [Gi]).

\begin{thm} Let $\Lambda_i, i=1,2$ be two transversal Lagrangian
submanifolds in a symplectic manifold $(M,\omega), dim\, M=2n$ (i.e.
they are transversal at all intersection points).
Then $Ext^i_{\cal C_M}(W_{\Lambda_1},W_{\Lambda_2})=0$ for all
$i\ne n$, and 
$Ext^n_{\cal C_M}(W_{\Lambda_1},W_{\Lambda_2})={\bf C}^{|\Lambda_1
\cap \Lambda_2|}$.

\end{thm}

Let us  consider the simplest case
of $1$-dimensional Lagrangian subspaces in the standard symplectic
${\bf R}^2$. 
We will be considering
algebraic quantization, not a smooth one.
Let $p,q$ be coordinates in ${\bf R}^2$ such that
$\{p,q\}=1$, and $A$ be the Weyl algebra (standard quantization of
this symplectic manifold, so that $[p,q]=\hbar\cdot 1$).
Let $\Lambda_1$ be the line $q=0$, and $\Lambda_2$ be the line
$p=0$. Then $W_{\Lambda_1}=A/Aq$, and $W_{\Lambda_2}=A/Ap$.
Clearly we have the following free resolution of $W_{\Lambda_1}$:

$$E^{\ast}: A\to A\to W_{\Lambda_1}\to 0,$$
where the first map is a multiplication by $q$, and the second map is
the natural projection. Then $Hom(E^{\ast},W_{\Lambda_2})$ gives
the complex $W_{\Lambda_2}\to W_{\Lambda_2}$, where the only map is
given by the multiplication by $q$. It is clear that there is
an isomorphism of ${\bf C}[q]$-modules: $W_{\Lambda_2}\simeq {\bf C}[q]$.
Therefore the complex  $W_{\Lambda_2}\to W_{\Lambda_2}$
has trivial cohomology  $H^0$, and $H^1$
is isomorphic to  ${\bf C}[q]/q{\bf C}[q]\simeq {\bf C}$.
Hence the only non-trivial Ext-group is 
$Ext^1(W_{\Lambda_1},W_{\Lambda_2})={\bf C}$. This proves
the theorem in the simplest case. 

Similarly one can check that for a given $\Lambda$ there is an isomorphism

$\oplus_{i\ge 0}Ext^i(W_{\Lambda},W_{\Lambda})\simeq H_{DR}^{\ast}(\Lambda)$. 
More generally, let $W_{\Lambda_i,L_i},i=1,2$ be small modules,
corresponding to closed Lagrangian transversal submanifolds $\Lambda_i$ which
carry local systems $L_i$. Assume that the Lagrangian submanifolds
intersect transversally. Then 

$Ext^n(W_{\Lambda_1,L_1},W_{\Lambda_2,L_2})=
\bigoplus_{x\in \Lambda_1\cap \Lambda_2}Hom((L_1)_x,(L_2)_x)$,
and all other $Ext$-groups are trivial. 

In the case of 2-dimensional tori there are two different interpretations
of the category of modules over the corresponding quantum torus:

a) as the category of small modules;

b) as the category of modules over a quantized algebra.

To a Lagrangian submanifold with a local system on it, one
assigns (in both cases a) and b)) an object of the corresponding
category of modules. It is natural to expect that these
objects correspond to each other under the  equivalence
of categories a) and b).

Although these observations
 show certain similarity of the Fukaya category $F(M)$
with the category ${\cal C}_M$ of modules
over the quantized function algebra, there is a difference in their
structures. Namely, Maslov index is not visible in ${\cal C}_M$, and
there is no non-trivial $A_{\infty}$-structure on the latter category.
On the other hand, Maslov index appears
in the definition of the structure of  $A_{\infty}$-category on $F(M)$.

\begin{que} Is there an $\A$-extension of the category ${\cal C}_M$ which
involves ``graded" objects 
(similarly to graded Lagrangian manifolds
in the construction of $F(M)$)?

\end{que}

This question is interesting even in the linear case (i.e. when
$M$ is the standard symplectic ${\bf R}^{2n}$). If the answer to the question
is affirmative, then for any three Lagrangian
submanifolds intersecting transversally
 there is a generalized Yoneda composition
$Ext^{\ast}(W_{\Lambda_1},W_{\Lambda_2})\times 
Ext^{\ast}(W_{\Lambda_2},W_{\Lambda_3})\to 
Ext^{\ast}(W_{\Lambda_1},W_{\Lambda_3})$ which involves
``instanton corrections", i.e. counting of holomorphic polygons
with the edges belonging to ${\Lambda_i},i=1,2,3$.

\subsection{Fukaya category for Lagrangian foliations and Moyal formula}

Fukaya suggested in [Fu1] to construct
a version of the category $F(T^{2n})$  for
Lagrangian foliations. 
To a foliation he assigned the algebra of foliation,
to a pair of ``good''  foliations he assigned a bimodule over the
corresponding algebras. This bimodule is a kind of a Floer complex.
The bimodule is not projective over either of the algebras.

Fukaya proposed an $A_{\infty}$-structure on the 
 category of Lagrangian foliations. Objects of the category
 are Lagrangian foliations together with transversal
measures. Spaces $Hom(F_1,F_2)$
 are analogs of Floer complexes
contsructed for ``good'' foliations $F_i, i=1,2$.
 The structure of $A_{\infty}$-category is given by
the ``higher compositions'' $m_k, k\ge 1$. 
Composition map $m_2$ is defined by the formula
which is  similar to the formula for the Moyal $*$-product
on a symplectic torus (see for example [We]):

$$ (f*g)(x)=({1\over \pi \hbar})^{2n}\int_{{\bf R}^{2n}\times {\bf R}^{2n}}
e^{-2\pi i\omega (x-y,x-z)/\hbar}f(y)g(z)dydz.
$$

 The formula contains the summation of the 
exponents of symplectic areas of triangles, which
makes it similar to the formula for $m_2$.
 The latter contains the exponents of the symplectic areas
of {\it holomorphic} polygons. Let us recall Fukaya's 
formula in the case when Lagrangian foliations are obtained from real
vector subspaces in the standard symplectic vector space ${\bf R}^{2n}$.

If one has three affine Lagrangian foliations: $F_i, i=1,2,3$, such that the
Maslov index of this triple is zero, then the formula looks
such as follows:

$$ m_2(\tau_2)(f\otimes g)(x,c,z)=\int_{F_2}\int_{F_2}
e^{-Q(a,b,c,\omega)}f(x,a,y)g(y,b,z)d\tau_2(F_2).
$$

The notation in the formula is explained in [Fu1]. Roughly speaking, 
$\tau_2$ is a transversal measure for $F_2$, so it
determines a section of $|\Lambda^{top}T(T^{2n})/TF_2|$.
Functions $f$ and $g$ live on the holonomy groupoids of
the pairs of foliations. The most interesting datum is the
function $Q$ which is defined as

$$
Q(a,b,c,\omega)=\int_{\Delta_{a,b,c}}\omega,
$$

where $\Delta_{a,b,c}$ is the geodesic triangle with vertices
in $a,b,c \in {\bf C}^n$.
This leads to the multi-dimensional theta-functions (see [Fu1]).

\begin{que} Can one explain the similarity of these two formulas
from the point of view of degeneration of the Fukaya category?
\end{que}

Apparently, the Fukaya category and the category of certain modules
over a quantized algebra of functions belong to the same connected
component of the ``moduli space of $A_{\infty}$-categories".

\section{The ``large complex structure" limit}

The content of this subsection is borrowed from  [KoSo2].

Non-commutativity can appear as a result of degeneration of 
the complex structure into a foliation.
The conventional approach to mirror symmetry suggests to consider the ``large complex
structure" limit, so that conjecturally  Calabi-Yau manifolds become foliated
by special Lagrangian tori (see [SYZ]). This gives rise to a 
stratum in the compactified ``moduli space of conformal
field theories. It is described
in terms of ``classical" theories.
The formal neighborhood
 of a classical theory in the moduli space 
 can be reconstructed from classical data by means of some formal rules.

 In this subsection we will discuss an
 example of this kind. 
For more details and applications to mirror symmetry see [KoSo2].
In  [GW] the conjecture similar to ours was formulated independently.
It was verified in [GW] in the case of K3 surfaces.

Here is the idea. Suppose that we study degenerations of a family of
$n$-dimensional Calabi-Yau manifolds $X_{\varepsilon}$ as $\varepsilon \to 0$.
Every $X_{\varepsilon}$ carries the Calabi-Yau metric 
$g(\varepsilon)=g_{ij}(\varepsilon), 1\le i,j\le N$. Let us rescale 
$g(\varepsilon)$ in such a way that the diameter of $X_{\varepsilon}$
will be $O(1)$ as $\varepsilon \to 0$.
We  consider the limiting space $X_0$, 
 where the limit is taken in Gromov-Haudorff metric.

The following conjectural description of $X_0$ was suggected in [KoSo2]:

1) $X_0$ is a metric space. It contains 
an  $n$-dimensional Riemannian manifold $X_0^{sm}$.
The dimension of $X_0^{sing}=X_0\setminus X_0^{sm}$ is less or equal
than two.

2) The manifold $ X_0^{sm}$ carries an
affine structure (i.e.  flat torsion-free 
connection on the tangent bundle $T_{X_0^{sm}}$).

3) There is a covariantly flat lattice $\Gamma\subset T_{X_0^{sm}}$.
In affine local coordinates it is given by ${\bf Z}^n\subset {\bf R}^n$.

4) Let us identify locally $ X_0^{sm}$ with ${\bf R}^n$. Then the Riemannian
metric $g$ on $ X_0^{sm}$ is K\"ahler-Einstein. This means that
 $g=\partial ^2H$ for some function $H$, such that
 $det(g)=const$ (the Monge-Ampere equation).
 
These data give rise to a fibration of flat tori on $X_0^{sm}$.

This conjectural picture is related to the mirror symmetry in
the following way. Consider (locally) the graph $dH$ in
$V={\bf R}^n\oplus {\bf R}^{\ast n}$. The latter space is a symplectic
manifold equipped with two Lagrangian foliations arising from
the coordinate spaces.
The graph $dH$ is a Lagrangian submanifold
in it. Since $H$ is defined up to the adding of an affine function,
the graph itself is defined up to translations.  Translated submanifolds are
still Lagrangian in $V$, so we will not pay much attention to this ambiguity.
Let $p_i,i=1,2$ be the canonical projections of $V$ to the coordinate subspaces.
Then the Monge-Ampere equation corresponds to the condition
$p_1^{\ast}(vol_{{\bf R}^n})=p_2^{\ast}(vol_{{\bf R}^{\ast n}})$
where $vol$ denotes the standard volume form.

Now we see that the whole picture is symmetric, so we can interchange
dual affine structures. Then $H$ gets replaced by its Legendre 
transformation. It does not change $X_0^{sm}$, and the limiting metric $g$.
It interchanges the dual affine structures and the dual lattices.
Hence it interchanges the corresponding dual fibrations of the flat tori.
This duality is geometric mirror symmetry. Conjecturally,
 degenerations of families of dual Calabi-Yau manifolds 
in the limit of ``large complex structure", lead to the two dual
fibrations of flat tori over the same base $X_0^{sm}$.

From the point of view
of $N=2$ superconformal field theory,
 the whole picture is classical, not quantum. 
 It is  shown in [KoSo2] how it simplifies
 the counting. In turns out that the counting of pseudo-holomorphic
 discs can be reduced to the counting of  certain binary trees in  $X_0^{sm}$.
It is also explained in [KoSo2] that the geometric degeneration is compatible
with certain degenerations of $N=2$ superconformal field theories.

 \begin{exa}
 
Moduli space of conformal field theories with the central charge $c=2$
 is a product of two modular
curves. These theories can be described as sigma-models
with a target space, which is a 2-dimensional torus.
Two separate ``infinite" limits for each modular curve give
rise to the two moduli spaces. 
 Maximal degenerations give rise to bundles
over a circle with the fibers which are circles themselves.
The mirror symmetry relates two dual circle bundles over the same base.
\end{exa}

\subsection{Non-commutative stratum}

As we have seen above, degenerations of the complex structure
can be described in terms of either commutative or non-commutative geometry.
It would be interesting to understand what kind of non-commutative theories
can be obtained in this way.

 Let us consider the case
of $K3$ surfaces  elliptically fibered
over ${\bf P}^1$. We can deform the complex
structure on a surface in such a way, that it becomes a foliation
on each non-degenerate fiber, and singular foliation on the degenerate
fibers and on the base (which is a $2$-sphere $S^2$). 
It is expected that to such geometric picture one can assign
a family of ``generically non-commutative" theories on
quantum tori (foliated fibers). 

\begin{que} How to describe these theories?

\end{que}

One also expects that similar non-commutative degenerations exist for arbitrary
Calabi-Yau manifolds. Loosely speaking, one can make into
non-commutative all special Lagrangian
tori in [SYZ] picture, thus adding a new non-commutative stratum to the moduli
space of $N=2$ superconformal theories.

\section{Other related topics}

\subsection{Quantum theta functions}

Let $T(L,exp(2\pi i\varphi))$ be a quantum torus, $\varphi\in {\bf R}$.
We fix a quadratic form $Q:L\times L\to {\bf C}$
with negative imaginary part
and linear functional $l:L\to {\bf R}$.
Let us also fix a symmetric bilinear form (,) on $L_{\bf R}$.
 We write $Q(x)=
(\Omega x,x)$ for some matrix $\Omega$ from the Siegel upper-half
space. In particular $\Omega$ is symmetric with respect to
the bilinear form.

 To these data we
associate a {\it quantum theta function} (see [M1]):

$$
\theta (Q,l,\varphi)=\sum_{a\in L}e^{2\pi i(Q(a)+l(a))}e(a).
$$

For any $t\in Hom(L, {\bf C}^{*})$ we define an automorphism
of $T(L,2\pi i{\varphi})$ by the formula:

$$
t^*(e(a))=t(a)e(a).
$$

Every element $\xi \in L$
defines $t_{\xi}\in Aut(T(L,2\pi i{\varphi}))$ such that
$t_{\xi}(e(a))=exp(2\pi i(-(a,\Omega \xi)+\varphi(a,\xi)))e(a)$.

The  following result is straightforward (see [M1]).

\begin{prp}
For an arbitrary $\xi$ as above we have:

$$t_{\xi}^*(\theta)=e^{-2\pi i(Q(\xi)-l(\xi))}e(\xi)\theta.$$

\end{prp}

\begin{que}  Let us write $Q(a)=(\Omega a,a)$, $\varphi(a,b)=
(\Phi(a),b)$. Is there an analog of the functional equation for
$\theta=\theta(\Omega,l,\Phi)$?

\end{que}
 We plan to return to this question
in the next paper.
Notice that standard proofs of the functional equation 
for theta-functions which use the Poisson summation 
formula do not work in quantum case.

\subsection{Quasi-modular forms}
 
There is an approach to the mirror symmetry for elliptic curves
due to Dijkgraaf (see [Di1]). It is the ``elliptic" version
of the counting of higher genus curves in the Calabi-Yau manifold.
 Main result of [Di1] is a statement that the generating function
counting certain coverings of the elliptic curve is a 
quasi-modular form in the sense of Kaneko and Zagier (see [KZ]). 

Dijkgraaf pointed out that for a fixed modular parameter $\tau$ and fixed 
K\"ahler class $t$ the related quantum field theory is based
on a four-dimensional  lattice of signature $(2,2)$ in ${\bf C}^2$
 (see [DVV]). 
Thus one has the group $O(2,2,{\bf Z})$ as the duality group 
of the quantum theory.
This group contains ${\bf Z}_2$
as a subgroup. The latter is responcible for the mirror symmetry.
The theory depends on a flat metric on a $2$-dimensional torus
and a skew-symmetric bilinear form on the lattice defining the torus.

The partition function of the theory is 
 quasi-modular. It can be modified,
so that it enjoys modular properties, but becomes non-holomorphic.

More precisely, for any $g>1$ Dijkgraaf defines a formal power series

$$ F_g(q)=\sum_{d\ge 1}N_{g,d}q^d,$$

where $q=exp(2\pi i \rho)$.
and $N_{g,d}$ counts the virtual number of ramified coverings of a complex
elliptic curve ${\cal E}$ by smooth complex curves of genus $g$. There are
$2g-2$ ramification points of index $2$, and coverings have degree $d$.
One can interpret $\rho$ as the complexified area of the underlying torus.

The following result can be found in the cited paper by Dijkgraaf.
It was rigorously proved in [KZ].

\begin{prp}
 $F_g\in {\bf Q}[E_2,E_4,E_6]$ and has weight $6g-6$,
where $E_k$ are Eisenstein
series.

\end{prp}

 Since $E_2$ is not a modular function, same is true for
$F_g$ . But $E_2$ is quasi-modular in the sense of [KZ], so does $F_g$.

It is natural to ask about the
non-commutative version of the Proposition.
If it exists, then the corresponding counting function should be
be a kind of  non-commutative limit of $F_g$. 
It is interesting to understand whether quantum tori
can be related to certain degenerations of quasi-modular forms.

Informally,  ``boundary modular forms''
 should correspond to the limits of quasi-modular forms
as modular parameter approaches to a  real number. It is natural to ask how
this limit should be understood. For example, one can consider limiting
functions as distributions or hyperfunctions.

The question about such a limit might be related to a different
question about ``boundary" modular forms in the sense of Zagier.
In Zagier's work, rational points of the boundary line
of the upper-half plane play a special role.
More precisely, he constructs embeddings of various spaces
of ``honest'' modular (or, hopefully, quasi-modular) forms to the space
of ``smooth function'' on ${\bf Q}$, modulo ``rational
functions''. The image of this embedding
consists of functions having modular properties.

Let us consider
an example of a ``boundary modular form" (I learned
it from Don Zagier).
Let $s(p,q)$ be a unique (after normalization)  function
defined for relatively prime integers $p,q$ in the following way:

a) $s(p,q) \in {\bf Q};$

b) $s(p+nq,q)=s(p,q);$

c) $s(-p,q)=-s(p,q);$

d) $s(p,q)+s(q,p)=(q^2+p^2+1-3pq)/12pq.$

Then as a function of $x=p/q$ it satiesfies the following
conditions:

1) $s(x+1)=s(x);$

2) $s(-1/x)=s(x)+P(x)$,

where $P(x)$ is the RHS of d). 

Then one can give a precise meaning to the
following statement:

the function $s(x)$ is modular as an element of the set of functions
on ${\bf P}^1({\bf Q})$, modulo  ``smooth'' rational functions.

\begin{que} What is the image of the Dijkgraaf's partition
function under the ``boundary embedding"?

\end{que}

From the point of view of the classical theory of
Eichler, Manin and Shimura,
the boundary modular forms might be related to
integrals of modular forms over geodesics between cusps in Lobachevsky
plane.
It is natural to ask about the meaning
of other geodesics.
One can expect that geodesics between arbitrary
boundary points should be treated within
the framework of the ``non-commutative'' geometry. 
Indeed, for non-cuspidal points
the geodesics are dense in the corresponding modular curve.

\subsection{On the Morita equivalence and deformation quantization}

Hopefully Morita equivalence of quantum tori can be generalized to 
other quantized function algebras.
It should be a statement that quantized algebras
corresponding  to different Poisson structures
and different values of the quantization parameter
produce  Morita equivalent algebras. 
Notice, however, that the
conventional deformation theory deals with formal series in 
a parameter $\hbar$.
From this point of view the transformation $\hbar\mapsto -1/\hbar$ does not
have sense. In order to allow such transformations (obviously we need them
in order to treat Morita equivalent quantum tori) we need to have a
``global'' moduli space of quantizations, not just a formal scheme
over $k[[\hbar]]$. The problem cannot be resolved at the level
of Poisson structures. It is ``quantum" problem. One can have
a discrete  group of symmetries of the quantum problem, which do not exists at
the level of Poisson structures.

\begin{que} Are there  examples of deformation quantization problems 
which model this phenomenon?
\end{que}

 In quantum groups
one often has series convergent  in $\hbar$. Let $q=exp(2\pi i \hbar)$.
It is interesting to look at the quantized homogeneous spaces 
(for example quantum groups themselves)
in the case $|q|=1$.

For example, let us consider the categories of projective modules
over the quantized coordinate rings ${\bf C}[G]_q$, where $q=e^{2\pi i/l}$,
where $l$ is a positive integer number,
and $G$ is a simple complex Lie group. One can see that
the quantized coordinate ring becomes a projective module
of finite rank $N=l^{dim(G)}$ over the subalgebra of the center,
which is isomorphic to the algebra ${\bf C}[G]$ of functions
on the Poisson-Lie group $G$.

\begin{que}
Is it true that for all primitive roots of $1$ of
the same order the quantized coordinate rings are Morita equivalent?

\end{que}

We would like to know the answer to a more general question. 
Namely, what is the ``duality group"
acting on the quantization parameter,
such that if parameters $q$ and $q^{\prime}$ 
belong to the same orbit, then the
corresponding algebras of functions
on the quantized homogeneous $G$-spaces 
${\bf C} [X]_q$ and ${\bf C} [X]_{q^{\prime}}$ are Morita equivalent?
More generally, the duality group should act on Poisson structures as well.

In the case of quantized coordinate rings of simple
Lie groups at roots of unity, the centers are non-isomorphic,
so one cannot expect the  group $SL(2,{\Z})$  be
the ``duality group'' of the theory. If the answer to the
above question is positive, then the ``duality group'' in question is the
 Galois group $Gal( {\bf Q}^{ab}/{\bf Q})$, where 
${\bf Q}^{ab}$ is the maximal abelian extension of the field
of rational numbers.

This is a toy-model for the global dualities in quantum physics.
Indeed, we have local quantizations given by perturbative
series, and we search for a global non-perturbative theory.

Another interesting question is to find an analog of the mirror symmetry
for the deformation quantization. 
In the case of quantum tori we discussed it in Section 4.

\vspace{10mm}

{\bf References}
 
\vspace{5mm}

[AP] D. Arinkin, A. Polishchuk, Fukaya category and Fourier transform,
math.AG/9811023.

\vspace{2mm}

[B1] S. Barannikov,  Semi-infinite Hodge structures 
and mirror symmetry for projective spaces, math.AG/0010157.

\vspace{2mm}

[B2] S.Barannikov, Quantum periods - I. 
Semi-infinite variations of Hodge structures, math.AG/0006193.

\vspace{2mm}

[B3] S. Barannikov, Generalized periods and mirror symmetry in dimensions $n>3$,
 math.AG/9903124.

\vspace{2mm}

[BK] S. Barannikov, M. Kontsevich,
Frobenius Manifolds and Formality of Lie Algebras of Polyvector Fields, 
alg-geom/9710032. 

\vspace{2mm}

[BEG] V. Baranovsky, S. Evens, V. Ginzburg
    Representations of quantum tori and double-affine Hecke algebras, 
   math.RT/0005024.
    
\vspace{2mm}

[BG] V. Baranovsky, V. Ginzburg, Conjugact classes in loop groups
and $G$-bundles on elliptic curves, alg-geom/9607008

\vspace{2mm}

[BM] P. Bouwknegt, V. Mathai, D-branes, B-fields and twisted K-theory,
 hep-th/0002023 

\vspace{2mm}

[BO] A. Bondal, D. Orlov,Reconstruction of a variety from 
the derived category and groups of autoequivalences, alg-geom//97120029.

\vspace{2mm}

[Co] A. Connes, Non-commutative geometry. Academic Press, 1994.

\vspace{2mm}

[CDS] A. Connes, M. Douglas, A. Schwarz, Non-commutative geometry
and Matrix theory: compactification on tori, hep-th/9711162.

\vspace{2mm}

[CK] Y-K. Cheu, M. Krogh, Non-commutative geometry from 0-branes in
a background $B$-field, hep-th/9803031.

\vspace{2mm}

[CR] A. Connes, M. Rieffel, Yang-Mills for non-commutative two-tori.
Contemporary Math. 62 (1987), 237-266.

\vspace{2mm}

[D] M. Douglas, D-branes and non-commutative geometry, hep-th/9711165.

\vspace{2mm}

[DFR] M. Douglas, B. Fiol, C. R\"omelsberger, Stability and BPS branes,
hep-th/0002037.

\vspace{2mm}

[Di1] R. Dijkgraaf, Mirror symmetry and elliptic curves. In: Moduli
Space of Curves. Progress in Mathematics, vol. 129, p.149-163.

\vspace{2mm}

[Di2] R. Dijkgraaf, Les Houches lectures on strings, fields
and duality, hep-th/9703136.

\vspace{2mm}

[DVV] R. Dijkgraaf, E. Verlinde, H. Verlinde, On moduli spaces
of conformal field theories with $c\ge 1$. In: Perspectives
in String Theory, P.Di Vecchia and J. L. Petersen Eds. World Scientific,
1988.

\vspace{2mm} 

[Fad] L. Faddeev, Modular double of quantum group, math.QA/9912078.

\vspace{2mm}

[FFFS] G. Felder, J. Fr\"ohlich, J. Fuchs, C. Schweigert,
Conformal boundary conditions and three-dimensional topological field theory,
 hep-th/9909140.
 
 \vspace{2mm}
 
[Fu 1] K. Fukaya, Floer homology of Lagrangian foliation and 
non-commutative mirror symmetry. Preprint 98-08, Kyoto University, 1998.

\vspace{2mm}

[Fu 2] K. Fukaya, Mirror symmetry of abelian varieties and multi theta
functions, preprint Kyoto University, 1998.

\vspace{2mm}

[FW] D. Freed, E. Witten, Anomalies in String Theory with D-Branes,
 hep-th/9907189.

\vspace{2mm}

[Gi] V. Ginzburg, Characteristic varieties and vanishing cycles,
Inv. Math. 84 (1986), 327-402.

\vspace{2mm}

 [GLO] V. Golyshev, V. Lunts, D. Orlov, Mirror symmetry for abelian varieties,
 math.AG/9812003.
 
 \vspace{2mm} 
  
[GW] M. Gross, P. Wilson, Large Complex Structure Limits of K3 Surfaces,
math.DG/0008018.

\vspace{2mm}

[Kap]  A. Kapustin, D-branes in a topologically nontrivial B-field,
  hep-th/9909089.

\vspace{2mm}

[KKO] A. Kapustin, A. Kuznetsov, D. Orlov,  Non-commutative instantons
and twistor transform, hep-th/0002193.

[Ko1] M. Kontsevich, Homological algebra of Mirror symmetry. Proceedings
of the ICM in Zurich, vol. 1, p. 120-139.

\vspace{2mm}

[Ko2] M. Kontsevich, Deformation quantization of Poisson manifolds, I, 
q-alg/9709040.

\vspace{2mm}

[Ko3] M. Kontsevich, Operads and Motives in Deformation Quantization,
 math.QA/9904055.
 
 \vspace{2mm}

[KoSo1] M. Kontsevich, Y. Soibelman, Deformation theory (book in preparation).

\vspace{2mm}

[KoSo2] M. Kontsevich, Y. Soibelman, Homological mirror symmetry and
torus fibrations, math.SG/0011041.

\vspace{2mm}

[KoSo3] M. Kontsevich, Y. Soibelman, Deformations of algebras over operads
and Deligne conjecture, Math. Phys. Stud., v. 21, 255-307 (2000).

[KR] M. Kontsevich, A. Rosenberg, Non-commutative smooth spaces, 
math.AG/9812158.

\vspace{2mm}

[KS] A. Konechny, A. Schwarz, $1/4$-BPS states on non-commutative tori,
hep-th/9907008.

\vspace{2mm}

[KZ] M. Kaneko, D. Zagier, A generalized Jacobi theta function and
quasi-modular forms. In: Moduli
Space of Curves. Progress in Mathematics, vol. 129, p.165-172.

\vspace{2mm}

[KoS] L. Korogodskii, Y. Soibelman, Algebras of functions on quantum
groups. Part I. Mathematical Surveys and Monographs, 56. American Mathematical
Society, Providence, RI, 1998. x+150 pp. ISBN: 0-8218-0336-0. 

\vspace{2mm}

[LLS] G. Landi, F. Lizzi, R. Szabo, String geometry and the non-commutative
torus, hep-th/9806099.

\vspace{2mm}

[M1] Yu. Manin, Quantized theta-function. Preprint RIMS, RIMS-700, 1990.

\vspace{2mm}

[M2] Yu. Manin, Mirror symmetry and quantization of abelian varieties,
math.AG/0005143.

[Mo] G. Moore, Finite in all directions, hep-th/9305139.

\vspace{2mm}

 [O1] D. Orlov,    On equivalences of derived categories of 
 coherent sheaves on abelian varieties, alg-geom/9712017. 

\vspace{2mm}

[O2] D. Orlov, Equivalences of derived categories and K3 surfaces, 
alg-geom/9606006.

\vspace{2mm}
 
[P1] A. Polishchuk, $A_{\infty}$-structures on an elliptic curve,
math.AG/0001048. 

\vspace{2mm}

[P2]  A. Polishchuk, Homological mirror symmetry with higher products,
 math.AG/9901025. 

\vspace{2mm}

[PZ] A. Polishchuk, E. Zaslow, Categorical mirror symmetry: 
the elliptic curve, math.AG/9801119.

\vspace{2mm}

 [RS] M. Rieffel, A. Schwarz,
 Morita equivalence of multidimensional noncommutative tori, math.QA/9803057. 
 
 \vspace{2mm}
 
[S1] A. Schwarz, Morita equivalence and duality, hep-th/9805034.

\vspace{2mm}

[Sab] C. Sabbah, Syst\`emes holonomes d'\'equations aux $q$-diff\'erences.
In D-modules and microlocal geometry, 1992, 126-147.

\vspace{2mm}

[Sch] V. Schomerus, D-branes and deformation quantization, hep-th/9903205.

\vspace{2mm}

[So] Y. Soibelman, Non-commutative compactifications 
and theta functions (in preparation).

\vspace{2mm} 

[SW] N. Seiberg, E. Witten, String theory and non-commutative geometry, 
hep-th/9908142.

\vspace{2mm}

[SYZ] A. Strominger, S-T. Yau, E. Zaslow, Mirror symmetry is $T$-duality,
hep-th/9606040.

\vspace{2mm}

[We] A. Weinstein, Classical theta functions and quantum tori.
Publ. RIMS. Kyoto Univ., vol. 30(1994), p. 327-333.

\vspace{3mm}

Address: Department of Mathematics KSU,

Manhattan, KS 66506, USA

e-mail: soibel@math.ksu.edu

\end{document}